\documentclass[10pt,oneside,a4paper]{article}

\usepackage{amssymb}
\usepackage{amsmath}
\usepackage{amsthm}
\usepackage{color}
\usepackage{graphicx}
\usepackage{wrapfig}
\usepackage{float}
\usepackage{subfigure}
\usepackage{rotating}
\usepackage{cite}
\usepackage{caption}
\usepackage{hyperref}
\usepackage{enumitem}
\usepackage{varwidth}

\theoremstyle{remark}

\theoremstyle{definition}

\def\RR{\mathbb{R}}
\def\CC{\mathbb{C}}
\def\ZZ{\mathbb{Z}}
\def\NN{\mathbb{N}}

\def\Proj{\mathrm{Proj}}
\def\span{\mathrm{span}}

\def\WB{W\!B}

\newcommand{\comment}[1]{}
\let\oldmarginpar\marginpar
\renewcommand\marginpar[1]{\-\oldmarginpar[\raggedleft\footnotesize #1]
{\raggedright\footnotesize #1}}

\author{
Victor Linroth
\thanks
{
Department of Mathematics, Uppsala University, Box 480, 751 06 Uppsala,  (Sweden). {\tt victor.linroth@math.uu.se}.
}
}

\title{Breakdown of hyperbolicity for quasiperiodic attracting invariant circles in a family of three-dimensional Henon-like maps}

\begin{document}
	
\maketitle

\begin{abstract}
We numerically study quasiperiodic normally hyperbolic attracting invariant circles that appear for certain parameter values in a family of three-dimensional Henon-like maps. These parameter values make up contour segments in the parameter space where the invariant circles have constant irrational rotation number. At the edges of these segments we find a breakdown of the hyperbolicity of the invariant circle. We observe the collision and loss of smoothness of two of the invariant Lyapunov bundles while the Lyapunov exponents all remain distinct. This is very similar to the breakdown of quasiperiodic normally hyperbolic invariant circles studied in previous works that have mostly focused on skew product type systems along with a few other special types of systems.

The numerical tools we use for finding the invariant circles and calculating rotation numbers, Lyapunov exponents and bundles are based on the recently developed Weighted Birkhoff method. To apply all of these tools we need for the invariant circles to be attracting (or repelling) and for the system to be invertible. This is a severe restriction compared to alternative methods, but it is very numerically efficient and allows us to study even highly irregular circles.
\end{abstract}

\section{Introduction}

\noindent The study of invariant manifolds has been at the very heart of dynamical systems since it's inception, because the existence or nonexistence of them gives great insight into the long term behavior of a dynamical system. Along with the question of existence, researchers have also extensively studied the question of persistence, i.e. for what systems does the invariant manifold persist under small perturbations. Part of this research yielded the notion of a \textit{normally hyperbolic manifold}, which informally is an invariant $C^1$ submanifold where the growth rate of the transverse vectors dominate the growth rate of the tangent vectors. The key results from \cite{Hirsch_Pugh_Shub_70,Fen_72} is then that a compact normally hyperbolic manifold without boundary that is invariant under a diffeomorphism persists under small $C^1$ perturbations. Later it was shown in \cite{Mane_1978} that the converse is also true.

With the persistence of compact normally hyperbolic manifolds well understood, the question then becomes to identify the scenarios under which they can break down. Naturally this is a very broad subject so we restrict ourselves to the case of normally hyperbolic invariant circles, the simplest form of compact manifold. This case has been studied extensively in the context of quasiperiodic skew product systems, see e.g. \cite{Haro_Llave_06}. In this paper we will examine the breakdown of normally hyperbolic invariant circles in a non-skew product type system induced by a diffeomorphism. These types of system has not been studied as extensively as the skew product case but some research has been done, see e.g. \cite{Calleja_Figueras_12,Canadell_Haro_2017_Alg}.

When studying invariant circles we separate the cases where the internal dynamics has rational rotation number and irrational rotational number. The breakdown of normal hyperbolicity in the case of a rational rotation number is well understood. Rational rotation number and normal hyperbolicity implies the existence of an even number of periodic orbits connected by stable manifolds which together make up the invariant circle. Here the circle will break down exactly when some of the Lyapunov exponents of the transversal dynamics equals the exponent of the internal dynamics, see \cite{Fen_72}. The case of irrational rotation number is not is more subtle. Assuming that the invariant circle is $C^{1+\varepsilon}$ and then it follows by a classical result of Denjoy that the internal dynamics is topologically conjugate to a rigid rotation. Like in the periodic case it is possible for the circle to break down when some of the Lyapunov exponents of transversal dynamics equals the internal exponent which is $0$. But the quasiperiodic internal dynamics also introduce new possibilities for breakdown scenarios. In \cite{Haro_Llave_06,Figueras_Haro_15} one such scenario is identified. There they present some families of skew product systems where the transversal Lyapunov exponents remain bounded away from the internal exponent but hyperbolicity is lost due to transversal directions losing regularity and colliding with the tangent direction on measure zero set. A similar scenario has been observed in the conformally symplectic setting, see \cite{Calleja_Figueras_12}. These are very similar to the phenomenon observed in this paper.

In this paper we use the efficient computational tools developed in \cite{Das_17,Das_Yorke_18} to numerically study the breakdown of attracting invariant curves in non-skew product systems. Of course more general tools for computing normally hyperbolic invariant manifolds could be used, see e.g. \cite{Haro_16}, but the computational gain of the tools chosen show a great advantage in time and offer all the observables needed for characterizing the breakdown such as invariant bundles, Lyapunov exponents and regularity of invariant curve and so on.

\section{Theoretical background/Setting}

Throughout this paper we will be considering the family of quadratic maps $F_a:\RR^3 \to \RR^3$ defined by
\[
F_a\left(\begin{array}{c}x\\y\\z\end{array}\right) =
\left( \begin{array}{c}y\\z\\Bx + M_1 + M_2 y - z^2\end{array} \right),
\hspace{10mm} a = (B,M_1,M_2) \in \RR^3,
\]
which is one possible generalization of the classic H\'enon family to three dimensions. We note that the maps all have constant Jacobian $\det(DF_a) = B \neq 0$ and are invertible. This family has also been considered in \cite{Gon_05}, where they establish more properties of the maps. Throughout this paper we will assume $B$ to be fixed and $|B| < 1$, i.e. $F_a$ is a dissipative system.

We will be interested in identifying $F_a$-invariant circles, these can be identified as curves $\mathcal{C} \subseteq \RR^3$ parametrized by an injection $K:S^1 \to \RR^3$ such that $K(S^1) = \mathcal{C}$ and $F_a \circ K = K \circ f$. Since $F_a$ is invertible the induced dynamics $f:S^1 \to S^1$ satisfying $F \circ K = K \circ f$ is invertible so $f$ is a circle homeomorphism whose dynamics are well understood. We will assume $f$ is orientation preserving, since if it is not we can look at $f \circ f$. Then we know that $f$ has a rotation number $\rho$ defined by
\[
\rho = \lim_{n\to\infty} \frac{\tilde{f}^n(x) - x}{n},
\]
where $\tilde{f}$ is a lift of $f$, i.e. any map from $\RR$ to $\RR$ satisfying $f \circ \pi = \pi \circ \tilde{f}$, where $\pi: \RR \to \RR/\ZZ = S^1$ is the projection $\pi(x) = x + \ZZ$.

The question now is when does these invariant circles exist? In \cite{Gon_05} they identify a curve in the $(M_1,M_2)$ space where the map has multipliers $e^{\pm i \varphi}$, $0<\varphi<\pi$ with $2\cos(\varphi) = -(M_2+1)/B$. They indicate that the conditions for both subcritical and supercritical Neimark-Sacker bifurcations are fulfilled. In particular we are interested in the supercritical bifurcation that gives birth to an attracting invariant circle. Near the Neimark-Sacker bifurcation there will be Arnold tongues with rational rotation numbers corresponding to values of $\varphi$. We are interested in invariant circles with irrational rotational number, in particular invariant circles that are persistent along some contour in the subspace. These do indeed exist between the Arnold tongues, and some justifications for this can be found in \cite{Canadell_Haro_2017} which establishes existence of quasiperiodic normally hyperbolic invariant tori in smooth families of analytic systems. We say more about normal hyperbolicity further down. The proof uses methods from Kolmogorov-Arnold-Moser (KAM) theory and as such demands we restrict ourselves to rotation numbers that are Diophantine.

A real number $\rho$ is Diophantine of type $(D,\nu)$ where $D>0$ and $\nu>2$ if
\[
\left| \rho - \frac{m}{n} \right| > D \frac{1}{|n|^{\nu}} \hspace{4mm} \textrm{for all } m,n \in \ZZ, \hspace{2mm} n \neq 0.
\]
A number is Diophantine if it is Diophantine of some type $(D,\nu)$, and we note that the set of Diophantine numbers has full Lebesgue measure in $\RR$.

We consider only invariant circles with quasiperiodic inner dynamics that are attracting and normally hyperbolic. This means that there is a splitting $T_{\mathcal{C}}\RR^3 = N_{-} \oplus L$, such that $\pi_{N_{-}}:N_{-} \to \mathcal{C}$ and $\pi_L:L \to \mathcal{C}$ are invariant subbundles of $T_{\mathcal{C}}\RR^3$, $L = DK(TS^1)$ and there are numbers $0 < \lambda_{-} < 1$ and $C_{-} > 0$ such that for any $x \in \mathcal{C}$ and any vector $v \in (N_{-})_x = \pi_{N_{-}}^{-1}(x)$
\[
\|DF_a^n(x) v\| \leq C_{-} (\lambda_{-})^n \|v \| \hspace{4mm} \textrm{for all } n \geq 0.
\]
Here $\| \cdot \|$ is the regular Euclidian norm on $T_{\mathcal{C}}\RR^3$. For a general definition of normal hyperbolicity see \cite{Hirsch_Pugh_Shub_70}. If $\mathcal{C}$ is a normally hyperbolic and attracting invariant circle with quasiperiodic inner dynamics we say it is completely reducible if there exists everywhere nonzero sections $e_0 \in \Gamma(L)$, $e_1,e_2 \in \Gamma(N_{-})$ such that $F_a e_0 = e_0$, $F_a e_1 = \lambda_1 e_1$ and $F_a e_2 = \lambda_2 e_2$ for real $\lambda_1,\lambda_2$ such that $ 0 < |\lambda_2| < |\lambda_1| < 1$.

Lastly we note that using a method referred to as \textit{bootstrap of regularity} it can be proved that if $\mathcal{C}$ is normally hyperbolic and the parametrization $K$ is $C^0$ then $K$ also has the same regularity as the system, in our case analytic. For more details see \cite[Section~3.2]{Haro_Llave_06_RR}.

\section{Numerical methods}

\noindent There are several ways to approach the problem of studying invariant circles in a smooth dynamical system and a complete list of all methods is probably yet to be compiled. Some of the history calculating invariant tori and many of the important developments can be found in \cite[Chapter~5.1]{Haro_16}. One of the possible methods to be considered is the parametrization method. The parametrization method is a very robust method for computing high order approximations and can be used for more general invariant manifolds than just tori, however we chose to forego the parametrization method for a couple of reasons. One reason is that the circles we seek to study are always attracting and quasiperiodic and so it is very simple and computationally efficient to approximate the torus by iteration of a point in its basin of attraction. We will informally refer to this as the 'iteration method' from now on. Another reason for not using the parametrization method is that it can be sensitive to irregularity of the invariant circle. This is because the method is based on Fourier expansion of the circle, and if the torus becomes highly irregular for some parameter values it can require a very large number of terms in the Fourier expansion to achieve satisfactory accuracy. The iteration method also suffers from a highly irregular circle. As we will see later the Lyapunov exponents associated to the attracting directions stay bounded away from zero but these directions can approach the invariant direction of the internal dynamics as parameters change. As explained in \cite{Haro_Llave_06} the attracting dynamics remains exponential but the constant terms grow when the directions approach. So while a point will converge slower toward an irregular torus it is still an exponential convergence and only adds a bit of computational time and no extra requirements on storage.

With the iteration method we can thus achieve a satisfactory approximation of an orbit on the invariant circle. It then remains to find methods for computing desired quantities for the circle. In our case this will mainly be the rotation number of the internal dynamics, and the Lyapunov exponents and their associated invariant sections on the circle. The invariant sections can be computed by iteration of tangent vectors in a similar manner of the iteration method itself, and has exponential convergence for the same reason. The rest of the quantities can be reformulated as integral problems where we integrate on the torus with respect to its invariant measure. For this we will use a recently developed method referred to as the Weighted Birkhoff method. 

\subsection{Calculating rotation numbers}

Let us begin with an explanation of how to compute the rotation number of an invariant circle given only the points of an orbit on the circle.We base our approach on the methods described in \cite{Das_16}. First choose an appropriate plane to project the orbit to, along with a coordinate system in that plane. For simplicity of presentation we can view these coordinates as real and imaginary parts and the plane as the complex plane $\CC$. Let $\{\gamma_n\}_n$ be the projected orbit and let $\gamma:S^1 \to \CC$ be a map such that $\gamma_n = \gamma(n\rho)$, where $\rho \in [0,1)$ is the rotation number we are looking to compute. Such a map will exist as long as the original dynamical system is of high enough differentiability.

Now pick a point $P \in \CC$ such that $P \notin \gamma(S^1)$, and choose a continuous function $\varphi:\RR \to \RR$ such that
\[
e^{2\pi i \varphi(x)} = \frac{\gamma \circ \pi(x) - P}{|\gamma \circ \pi(x)-P|}
\]
where $\pi:\RR \to S^1= \RR/\ZZ$ is the projection map. We note that $\varphi(x+1) = \varphi(x) + W$ where $W\in\ZZ$ is the winding number around $P$ and let $\tilde{\Delta}:\RR \to \RR$ be defined by $\tilde{\Delta}(x) = \varphi(x+\rho)-\varphi(x)$. Since $\tilde{\Delta}(x+1) = \varphi(x+\rho+1)-\varphi(x+1) = \varphi(x+\rho)+W-\varphi(x)-W = \tilde{\Delta}(x)$ we see there is $\Delta:S^1 \to \RR$ such that $\tilde{\Delta} = \Delta \circ \pi$. We know that the sequence $\{k\rho\}_{k=0}^{\infty}$ on $S^1$ is equidistributed with respect to Lebesgue measure $m$ on $S^1$, so with $\Delta_k = \Delta(k\rho)$, Birkhoff's ergodic theorem gives that
\[
\lim_{n\to\infty}\frac{1}{n}\sum_{k=0}^{n-1} \Delta_k = \int_{S^1} \Delta dm = \int_0^1 \tilde{\Delta}(x)dx = \int_0^1 (\varphi(x+\rho)-\varphi(x))dx
\]
\[
= \int_0^1 \int_0^{\rho} \varphi'(x+s)ds dx = \int_0^{\rho}\int_0^1\varphi'(x+s)dx ds = W\rho.
\]
In practice one usually tries to choose plane, projection and point $P$ in such a way that $W=1$. Computing the numbers $\{\Delta_k\}$ in practice requires a bit of work since we don't immediately have the function $\varphi$. To do this one first finds numbers $\{\Delta_k^*\}_{k=0}^N \subset [0,1)$ such that $\Delta_k^* =(2\pi)^{-1}( \arg(\gamma_{k+1})-\arg(\gamma_k))$  for some choice of $\arg$. Then $\Delta_k^*$  and $\Delta_k$ will only differ by a integers. We can modify $\{\Delta_k^*\}$ by integers such that if $k\rho$ and $m\rho$ are 'close' in $S^1$ then $\Delta_k^*$ and $\Delta_m^*$ should be close as well. What constitutes 'close' will depend on how the curve behaves, and for larger $N$ we can find closer points. If this process succeeds then $\{\Delta_k\}$ and $\{\Delta_k^*\}$ will differ by the same integer everywhere and we can find it by noting that we want $\rho \in [0,1)$.

One problem with determining when $k\rho$ and $m\rho$ are close in $S^1$ is that we do not have these points but only their images under $\gamma$. Since we are looking at a projection of an invariant circle it is possible that the projected circle is self intersecting, i.e. two different points map to the same point. The solution proposed in \cite{Das_16} is to use a delay coordinate embedding $\Gamma:S^1 \to \RR^{2L}$ defined by $\Gamma(x) := (\gamma(x),\gamma(x+\rho),\hdots,\gamma(x+(L-1)\rho)))$, for some $L \geq 2$. This gives a new orbit $\Gamma_n = (\gamma_n,\gamma_{n+1},\hdots,\gamma_{n+L-1})$ in $\RR^{2L}$, and when we want to determine if $k\rho$ and $m\rho$ are close we compare $\Gamma_k$ and $\Gamma_m$. The theoretical justification given in \cite{Das_16} for why this works is that by the Whitney and Takens Embedding theorems the map $\Gamma$ will be an embedding for almost every curve $\gamma$, so this mapping should be a reasonably effective method to remove self intersections.

With this we have reduced the problem of computing the rotation number to the problem of approximating the limit of the Birkhoff sums
\[
B_n^{\Delta} = \frac{1}{n}\sum_{k=0}^{n-1} \Delta_k.
\]
The problem now is that $\{B_n^{\Delta}\}$ generally only have $\mathcal{O}(1/n)$-convergence, which is unsatisfactory. A very elegant and simple solution to this is to instead look at weighted Birkhoff sums
\[
\WB_n^{\Delta} = \frac{1}{A_n}\sum_{k=1}^{n-1} w\left(\frac{k}{n}\right)\Delta_k, \hspace{6mm} A_n = \sum_{k=1}^{n-1} w\left(\frac{k}{n}\right),
\]
where the weight $w:R \to [0,\infty)$ is a smooth bump function with support $[0,1]$ and $\int_{\RR}w(x)dx \neq 0$. The particular weight function we use is
\[
w(t) = \left\{ \begin{array}{ll}
\exp\left(\frac{1}{t(1-t)}\right), & t \in [0,1] \\
0, & t \notin [0,1],
\end{array} \right.
\]
but any choice satisfying the conditions works. If $\Delta$ is a function and $\rho$ is Diophantine, then
\begin{equation*}
\left| \WB_n^{\Delta} - \int_{S^1} \Delta dm \right| = \mathcal{O}\left(\frac{1}{n^k} \right)
\end{equation*}
for any $k \in \NN$. In other words $\{\WB_n^{\Delta}\}$ converges super-polynomially. Here the constants in the bounds will depend on $k$, $\rho$, $\Delta$ and $w$. The proof for this is given in \cite{Das_Yorke_18}. We should note here that the constants for $\mathcal{O}\left(1/n^{k} \right)$ derived in the proof contains $\|w^{(k)}\|$ as a factor, and these grow very fast in $k$. Indeed they grow faster than $c^k k!$ for all $c>0$. We also note that in the case of non-Diophantine $\rho$, $\{\WB_n^{\Delta}\}$ still converges but only $\mathcal{O}(1/n)$.

Another thing worth mentioning is that in the particular case of calculating rotation numbers, the Weighted Birkhoff method also converges super-polynomially for rational rotation numbers, which can be of relevance for numerical application. We are not aware of this result being published anywhere else so we give an outline of the proof. Assume for simplicity that $W=1$ and let $\varphi(x) = x + \sum_{m \in \ZZ} a_me^{2\pi i m x}$. Then we see that
\[
\tilde{\Delta}(x) = \varphi(x+\rho) - \varphi(x) = \rho +  \sum_{m \in \ZZ} a_m(e^{2\pi i m \rho} - 1)e^{2\pi i m x}
\]
Since $\Delta_k = \Delta(k\rho) = \tilde{\Delta}(k\rho)$ we get that
\[
E_n := \frac{1}{A_n} \sum_{k=0}^{n-1} w\left(\frac{k}{n} \right) \Delta_k - \rho = \frac{1}{A_n} \sum_{k=0}^{n-1} w\left(\frac{k}{n} \right) \sum_{m \in \ZZ} a_m(e^{2\pi i m \rho} - 1)e^{2\pi i m k \rho}
\]
\[
= \frac{1}{A_n} \sum_{m \in \ZZ} a_m(e^{2\pi i m \rho} - 1) \sum_{k \in \ZZ} w\left(\frac{k}{n} \right) e^{2\pi i m k \rho}
\]
\[ = \frac{1}{A_n} \sum_{m \in \ZZ} a_m(e^{2\pi i m \rho} - 1) \sum_{k \in \ZZ} \int_{\RR} w\left(\frac{t}{n} \right) e^{2\pi i m t \rho} e^{-2\pi i k t} dt
\]
where in the last step we used the Poisson summation formula. If $\rho$ is rational then some of the terms are zero, so let $\ZZ^{(\rho)} := \left\{m \in \ZZ \, \big| \, m \rho \notin \ZZ \right\}$. Then doing a change of variables and integrating by parts $d$ times we get that
\[
E_n = \frac{n}{A_n} \sum_{m \in \ZZ^{(\rho)}} a_m(e^{2\pi i m \rho} - 1) \sum_{k \in \ZZ} \int_0^1 w(s) e^{2\pi i (m \rho - k) n s} ds
\]
\[
= \frac{n}{A_n} \sum_{m \in \ZZ^{(\rho)}} a_m(e^{2\pi i m \rho} - 1) \sum_{k \in \ZZ} \frac{1}{(2 \pi i (m\rho - k)n)^d} \int_0^1 w^{(d)}(s) e^{2\pi i (m \rho - k) n s} ds.
\]
Estimating all the terms we get that $E_n = \mathcal{O}(1/n^d)$, for any $d \geq 1$. For more details we refer to \cite{Das_Yorke_18} where they do the proof for more general functions and Diophantine rotation numbers. The reason it also works for our case with a rational rotation number is that our function has the special form $\varphi(x+\rho) - \varphi(x)$ which leads to cancelation of terms that otherwise would cause problems.

\subsection{Tracing contours in parameter space}
\label{subsec:tracing}

Once we have a reliable way of calculating rotation numbers we turn to the problem of finding the parameter values which give a certain rotation number. In other words if $R \subseteq \RR^2$ is the subset of the parameter space which gives systems with attracting invariant circles and we let $r:R \to [0,1)$ be the function which gives the rotation number, then we are looking for the contours $r^{-1}(\rho) = \{x \in R \mid r(x) =\rho\}$ for some given $\rho \in [0,1)$. In our case we know that the contours we are looking for will be $1$-parametric curves. If we can locate a candidate area for the search we choose a line segment $\ell:[0,1] \to R$ that transversally crosses the contours we are interested in. If done properly the function $r \circ \ell$ will be Devil's staircase type function and in particular it will be monotone. So by checking if $r \circ \ell(0) \leq \rho \leq r \circ \ell(1)$ or $r \circ \ell(1) \leq \rho \leq r \circ \ell(0)$ we know if the line segment crosses our desired contour.

Here there is a large variety of numerical methods for finding where a function on an interval attains a certain value, but we note that due to the Devil's staircase nature of our function we can not utilize any method based on taking derivatives. Instead we have chosen to use the False Position Method (FPM) which gives good performance for $\rho$ with good Diophantine properties since near such a number the 'steps' of the staircase are relatively small. The FPM is based on constructing a nested sequence of intervals $[a_k,b_k]$ enclosing the solution. So for $f(t) = r \circ \ell (t) - \rho$ and $a_0,b_0$ equal to $\ell(0), \ell(1)$ such that $a_0 \leq b_0$. Then given $[a_k,b_k]$ we compute
\[
c_k = \frac{a_k f(b_k) - b_k f(a_k)}{f(b_k) - f(a_k)}
\]
and depending on the sign of $f(c_k)$ we get either $[a_{k+1},b_{k+1}] = [a_k,c_k]$ or $[c_k,b_k]$. There are some problems when analyzing the convergence of the FPM, in particular there are examples of functions for which the method does not converge. However for our application we did not encounter any of these problems, and if one does encounter these problems there are modified versions of the method one can use.

With a good approximation $p_0 \in R$ of one point on the contour we turn to the problem of finding more. Since the contour continues in two directions we must decide on an initial direction to search. With this we start to examine the value of $r$ on points of a circle around $p_0$ with some prescribed small radius. Taking two distinct points on the circle in the chosen initial direction, we evaluate $r$ there and if these values enclose $\rho$ we are done. If not we start stepping along the circle until we find two points with values of $r$ that do enclose $\rho$. Which way we start stepping depends on the values of $r$ we got and which inclination of $r$ on the circle these values suggest. When we have an enclosure we take the corresponding points and look at the line segment between them. Here we can reuse the previously explained FPM to find another good approximation $p_1$ of a point on the contour.

From here we repeat the same procedure for the new points using the old points to guess the initial direction to search. Since the contour does not continue indefinitely the procedure will eventually fail, in which case we shrink the radius of the circle used to search in order to get good approximations near the point of breakdown.

\subsection{Lyapunov exponents and associated sections}

In order to calculate the Lyapunov exponents, we see that applying Oseledet's ergodic theorem to our (invertible) dynamical system with the derivative cocycle gives us a decomposition of the tangent spaces
\[
T_x\RR^3 = \RR^3 = \bigoplus_{i=0}^{s(x)} H_x^i \hspace{2mm} x \in \RR^3, \hspace{2mm} \textrm{ such that } \hspace{2mm} DF_a(x)H_x^i = H_{F_a(x)}^i, \hspace{4mm} i = 0,1,\hdots, s(x),
\]
and
\[
\lim_{n \to \pm\infty} \frac{1}{n} \log\| DF_a^n(x) v\| = \lambda_i(x), \hspace{4mm} v \in H_x^i\setminus\{0\}, \hspace{4mm} i = 0,1,\hdots, s(x),
\]
where $\lambda_0(x) > \hdots > \lambda_{s(x)}(x)$ are the Lyapunov multipliers. We will restrict ourselves to the basin of attraction of a normally hyperbolic attractor where the subspaces $\{H_x^i\}$ constitute an orientable subbundle $H^i$ of the Tangent bundle. In this case the Lyapunov exponents will not depend on the point $x$ and we will write $\lambda_0 > \hdots > \lambda_s$. Note that in the case of quasiperiodic attractor we get that $\lambda_0 = 0$ and we have either $s=1$ or $s=2$. To simplify exposition a bit we assume $s=2$, but the following methods works for the case $s=1$ too.

We are interested in numerically approximating $\lambda_i$ and $H_x^i$ for $x$ on the attractor, $i=0,1,2$. We begin by finding a procedure for approximating vectors spanning $\{H_x^i\}$. If we take a vector $u \in \RR$ and view it as a tangent vector for some point $x$ we see that $DF_a^n(x) u$ will be very close to $H_{F_a^n(x)}^0$ for large $n$, as long as $u$ was not in the subspace $H_x^1 \oplus H_x^2$ (and even if it was, it is possible for numerical inaccuracies to push it out of the subspace). With this in mind we look at the sequence of vectors defined by the Gramm-Schmidt scheme
\[
\left\{\begin{array}{l}
u_{k+1} = \|DF_a(x_k)u_k\|^{-1}DF_a(x_k)u_k, \\
v_{k+1} = \frac{DF_a(x_k)v_k - \Proj_{u_{k+1}}DF_a(x_k)v_k}{\left\|DF_a(x_k)v_k - \Proj_{u_{k+1}}DF_a(x_k)v_k\right\|}
\end{array} \right.
\]
where $x_{k+1} = F_a(x_k)$ is the orbit and $\Proj_{u}v = \|u\|^{-1}(u \cdot v) u$. We see then that $\span(u_k) \to H^0_{x_k}$ and $\span(u_k,v_k) \to H_{x_k}^0 \oplus H_{x_k}^1$ as $k \to \infty$. Since the dynamical system is invertible we can perform the same scheme backwards
\[
\left\{\begin{array}{l}
\tilde{u}_{k-1} = \|DF_a(x_{k-1})^{-1}\tilde{u}_k\|^{-1}DF_a(x_{k-1})^{-1}\tilde{u}_k, \\
\tilde{v}_{k-1} = \frac{DF_a(x_{k-1})^{-1}\tilde{v}_k - \Proj_{\tilde{u}_{k-1}}DF_a(x_{k-1})^{-1}\tilde{v}_k}{\left\|DF_a(x_{k-1})^{-1}\tilde{v}_k - \Proj_{\tilde{u}_{k-1}}DF_a(x_{k-1})^{-1}\tilde{v}_k\right\|}
\end{array} \right.
\]
where  $x_{k-1} = F_a^{-1}(x_k)$. Then we get that $\span(\tilde{u}_k) \to H_{x_k}^2$ and $\span(\tilde{u}_k,\tilde{v}_k) \to H_{x_k}^1 \oplus H_{x_k}^2$ as $k \to \infty$.

In practice we choose some appropriate $N_1,N_2,N_3 > 0$ and $N = N_1+N_2+N_3$ and some point $x_1$ that is approximately on the attractor, as well as some choice of vectors $u_1,v_1,\tilde{u}_N,\tilde{v}_N$. Then we compute the orbit $\{x_k\}_{k=1}^N$ and iterate the vectors $u_1,v_1$ forward starting at $x_1$ and the vectors $\tilde{u}_N,\tilde{v}_N$ backward starting at $x_N$. If $N_1$ and $N_3$ were large enough then the vectors $\{u_k,v_k,\tilde{u}_k,\tilde{v}_k\}_{k=N_1}^{N_1+N_2}$ should give us a good approximation of $H^0,H^1,H^2$ on the orbit segment $\{x_k\}_{k=N_1}^{N_1+N_2}$ and
\[
\left\{\begin{array}{l}
H_{x_k}^0 \approx \span(u_k), \\
H_{x_k}^1 \approx \span(u_k,v_k) \cap \span(\tilde{u}_k,\tilde{v}_k), \\
H_{x_k}^2 \approx \span(\tilde{u}_k).
\end{array} \right.
\]
For the Lyapunov exponents we note that if $w_1 \in H_{x_1}^i$ and $\|w_1\| = 1$, then we see that
\[
\lambda_i = \lim_{n \to \infty} \frac{1}{n} \log\| DF_a^n(x_1) w_1\| = \lim_{n \to \infty} \frac{1}{n} \sum_{k=1}^n \log\| DF_a(x_k) w_k\|
\]
where $w_{k+1} = \|DF_a(x_k) w_k\|^{-1}DF_a(x_k) w_k$ for $k \geq 1$. Since $H_x^i$ are one-dimensional and orientable we see that there is a section $w$ on the attractor such that $w(x_k) = w_k$, and $\|w(x)\| =1$ for all $x$. Then by Birkhoff's theorem we get that
\[
\lambda_i = \lim_{n \to \infty} \frac{1}{n} \sum_{k=1}^n \log\| DF_a(x_k) w(x_k)\| = \int_A \log\|DF_a(x)w(x)\| d\mu(x)
\]
where $A$ is the attractor and $\mu$ is the invariant measure induced by the dynamical system. When we have this type of integral over an invariant measure we see that we are in a setting where the Weighted Birkhoff method is applicable, and this is how we use our approximations of the sections spanning the invariant subspaces to approximate the Lyapunov exponents.

\section{Results}

\begin{figure}[htb!]
\centering
\includegraphics[width=.45\linewidth]{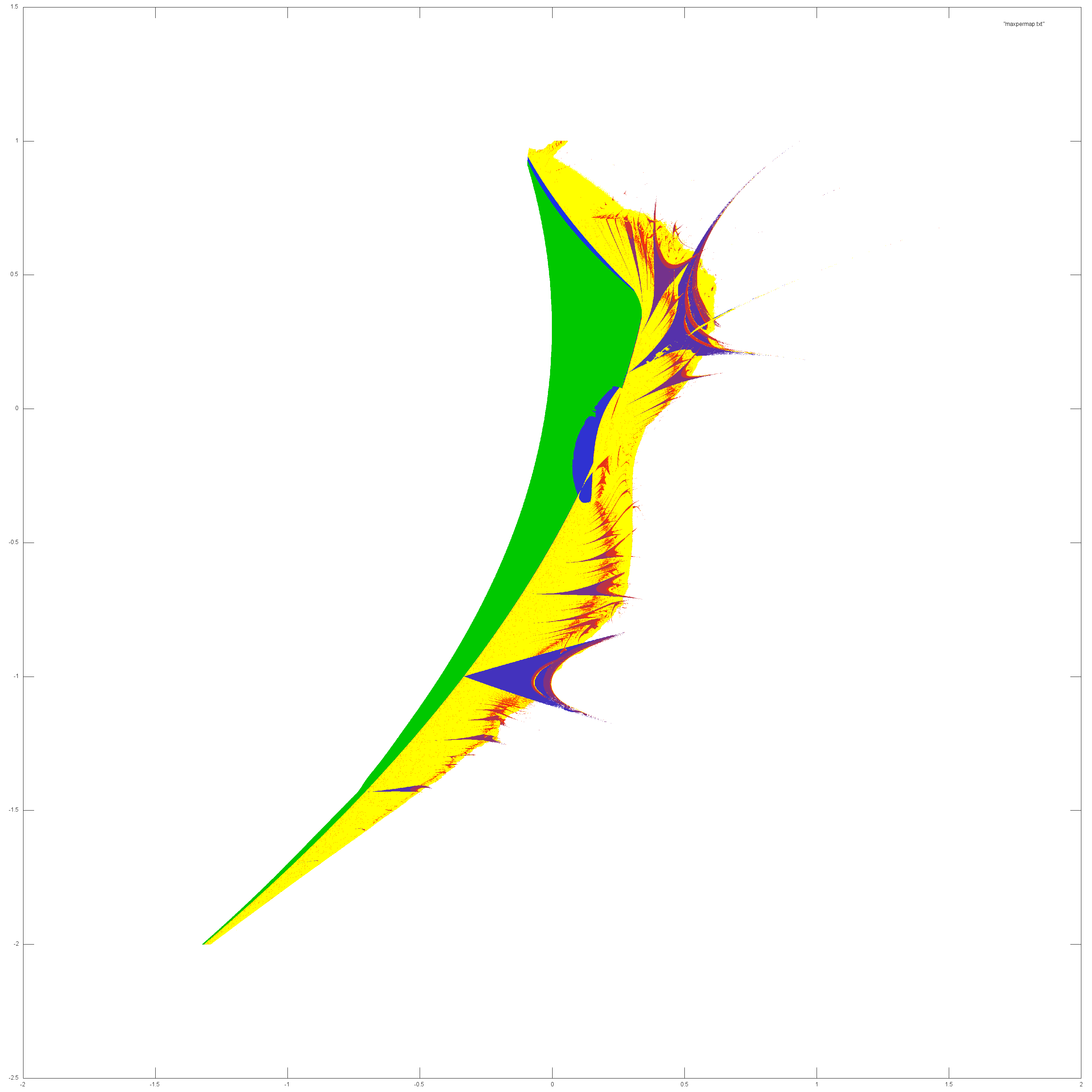}
\includegraphics[width=.45\linewidth]{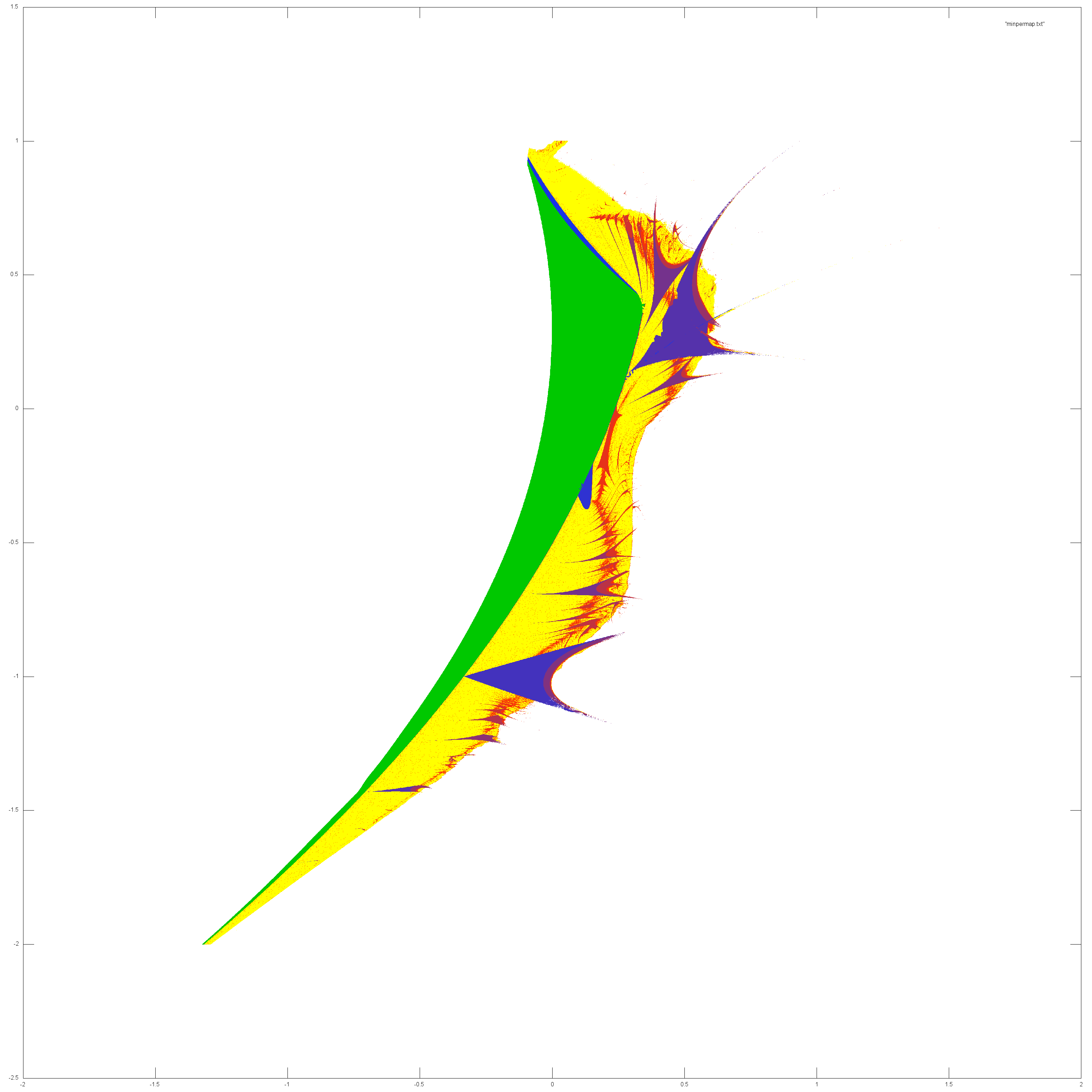}
\caption{\textit{Left:} Maximum period of attractor found for parameter value. \textit{Right:} Minimal period found. Color code: Green is fixed point. Blue is period 2 and higher periods goes from blue through purple to red. Yellow is no discernible period found and white is no attractor found at all.}
\label{fig:permap}
\end{figure}

The methods described in the previous section requires us to have a good idea of where the attractors we are looking for are, both in terms of the physical space that the dynamical system is operating on and i terms of the parameter space. Locating the attractor in physical space is not too hard and we can visually inspect it by plotting it for whichever parameter value we choose. To get a good idea of which parameters have which type of attractors we first choose some grid of values in the parameter space and for each value on the grid take a grid of points in the physical space. For each of these we iterate the point forward and test for periodic behavior. The results can be seen in Figure \ref{fig:permap}. First thing to note is that the picture is a bit different depending on if we plot the maximum period found or minimum period found. This indicates that there can be more than one attractor for each parameter value. We note that the pattern on the right that is very similar to the one found when studying the parameter space for the classic family of Arnold circle maps. Because of this we have good reason to believe there are curves of parameter values for which we have quasiperiodic attractors.

\begin{figure}
\center
\includegraphics[width=0.9\linewidth]{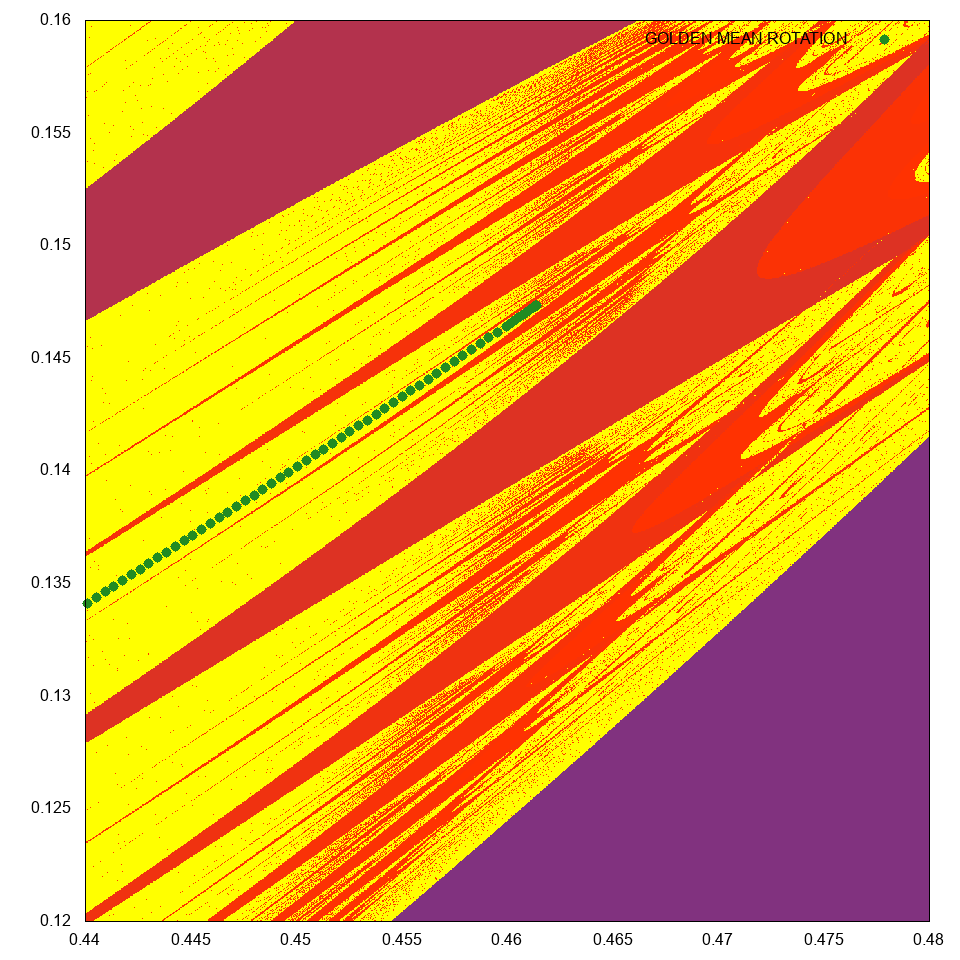}
\caption{Zoomed in picture of maximum period found. The Green dots are approximate points of the contour corresponding to attractors with rotation number equal to the golden mean.}
\label{fig:permapzoom}
\end{figure}

\begin{figure}
\center
\includegraphics[width=0.45\linewidth]{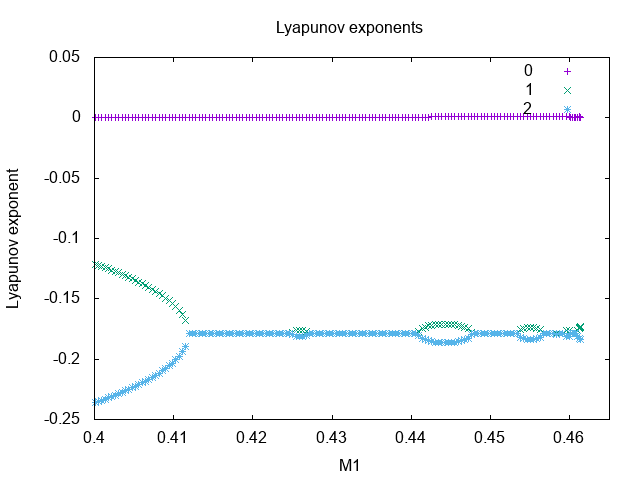}
\includegraphics[width=0.45\linewidth]{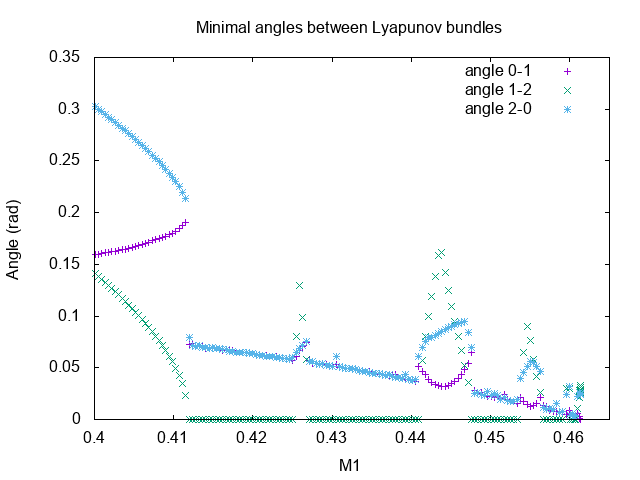}
\caption{\textit{Left:} Numerical approximations of Lyapunov exponents. Note that two of the exponents seem to coincide for some parameter values. \textit{Right:} Numerical approximations of minimum angles between Lyapunov bundles. In the cases where two exponents seem to coincide the corresponding angle has been set to 0, meaning they are the same (two dimensional) bundle.}
\label{fig:lyapunov}
\end{figure}

We begin with searching for attractors with rotation number equal to the golden mean $\varphi$, as this is the number with the ``best'' Diophantine properties and hence should give us good performance for our numerical methods. In fact the golden mean is what is called a badly approximable number, meaning there is some $D>0$ such that $\varphi$ satisfies
\[
\left| \varphi - \frac{m}{n} \right| > D \frac{1}{n^2} \quad \textrm{for all } m,n \in \ZZ, n \neq 0.
\]
Using the methods described in Section \ref{subsec:tracing} we locate a parameter point giving the desired rotation number and then trace the contour as it continues up to the right. These points are plotted in Figure \ref{fig:permapzoom} which is a zoomed in version of Figure \ref{fig:permap}.
As we are tracing the contour we are also calculating and saving the Lyapunov exponents as well as the minimum angles between the Lyapunov bundles of the attractors. These values are plotted in Figure \ref{fig:lyapunov}. We note that one of the Lyapunov exponents stay constant around zero which is expected since quasiperiodic motion has Lyapunov exponent zero, so this is the exponent for the tangential direction. The other one or two exponents are negative which we expect from an attractor, and we note that the stay away from zero. So if there is a loss of hyperbolicity at the end of the contour it is not due to Lyapunov exponents going to zero. We can also note that there seems to be three distinct exponents at the end, indicating that attractor is completely reducible at the breakdown. Looking at the minimal angles we can see that minimum between the tangential direction and the slow contracting direction is getting very small towards the end, while the other angles all stay away from zero.

\begin{figure}
\center
\includegraphics[width=0.45\linewidth]{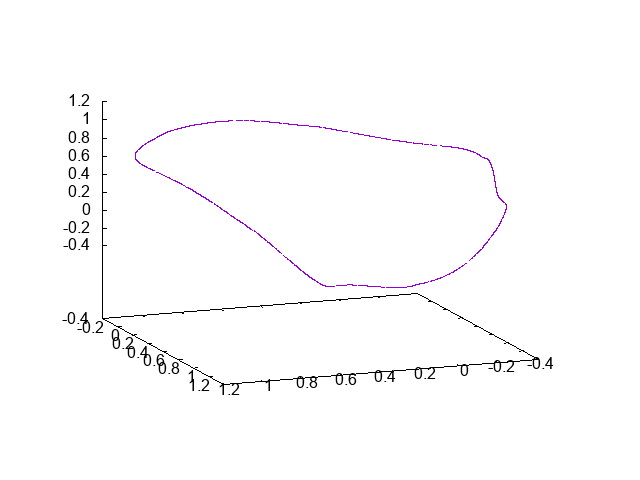}
\includegraphics[width=0.45\linewidth]{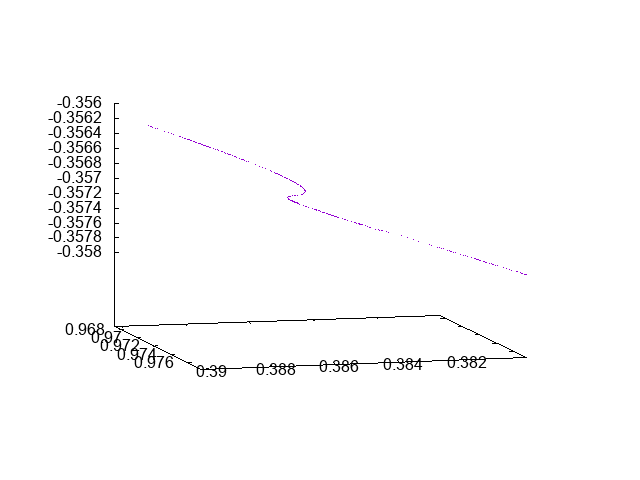}
\caption{\textit{Left:} Attractor at the last step. \textit{Right:} Zooming in we can see small 'folds' that are not visible in the right picture.}
\label{fig:attractor}
\end{figure}

\begin{figure}
\center
\includegraphics[width=0.45\linewidth]{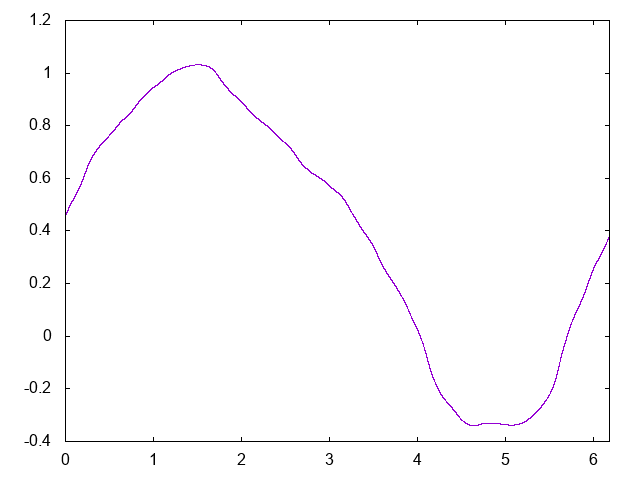}
\includegraphics[width=0.45\linewidth]{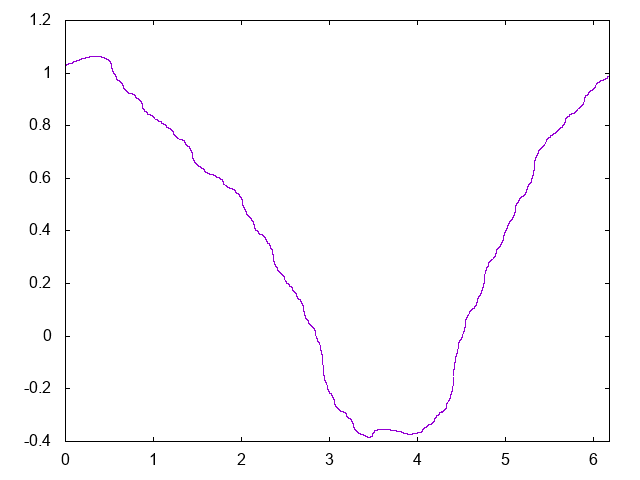}
\caption{\textit{Left:} Projection of the conjugation on one of the coordinate axes for parameter values some distance away from the breakdown. \textit{Right:} Projection at the last step.}
\label{fig:xcoord}
\end{figure}

\begin{figure}
\center
\includegraphics[width=0.45\linewidth]{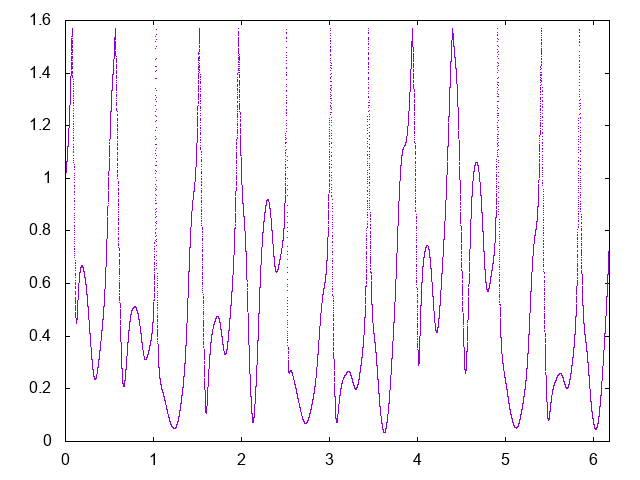}
\includegraphics[width=0.45\linewidth]{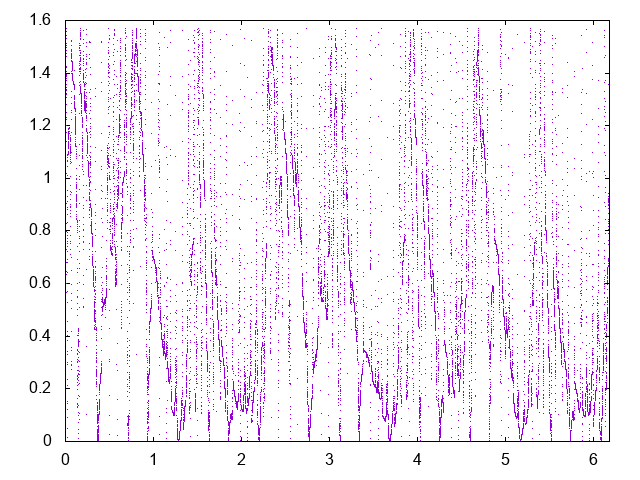}
\caption{\textit{Left:} Smallest angles between the tangent direction and the slow contracting direction for parameter values some distance away from the breakdown. \textit{Right:} Angles at the last step. }
\label{fig:minangle}
\end{figure}

In Figure \ref{fig:attractor} we see two plots of the attractor at the last step. When looking at the whole attractor it appears to have some amount of differentiability, although if we zoom in far enough we find small 'folds' along the curve so it is not quite as regular as it first appears. But even with these folds it does not appear the curve loses differentiability as the parameters approach the point of breakdown. Since we know the rotation number it is relatively simple to construct the conjugation. In Figure \ref{fig:xcoord} we show the projection of this conjugation on one of the coordinate axes. To the left we have the projection some distance away from the breakdown and it does not give us any good reason to believe it is anything less than analytic. To the right we have the projection at the last step and it appears as if differentiability is about to become lost. In Figure \ref{fig:minangle} we show the smallest angles between the tangent direction and the slow contracting direction and here we see indeed there appears to be a collision between the two.

We note here that if it is the case that the conjugation looses differentiability at the breakdown then this places distinct limits on how regular the curve itself can be. If the curve is $k>2$ times differentiable then the conjugation is $k-1-\varepsilon$ times differentiable for any $\varepsilon > 0$, c.f. \cite{katznelson_ornstein_1989,khanin_sinai_1987}. So if the conjugation is not differentiable then the curve can not be more than $2$ times differentiable. 

\section*{Acknowledgments.}
This work was supported by the Swedish Research Council Grant 621-2011-3629. I want to thank S. V. Gonchenko, \'A. Haro, J. Yorke and J.-L. Figueras for fruitful discussions. 

\bibliography{bibliography}{} 
\bibliographystyle{alpha}

\end{document}